\documentclass[11pt]{amsart}

\usepackage{amssymb,enumitem,url}

\usepackage[mathlines]{lineno}
\usepackage[colorlinks=false]{hyperref}

\usepackage{cleveref}

\theoremstyle{plain}

\newtheorem{theorem}{Theorem}[section]
\newtheorem{lemma}[theorem]{Lemma}
\newtheorem{proposition}[theorem]{Proposition}
\newtheorem{conjecture}[theorem]{Conjecture}
\newtheorem{sublemma}{}[theorem]

\newcommand{\mcal}[1]{\ensuremath{\mathcal{#1}}}
\newcommand{\ex}[1]{\ensuremath{\textup{EX}_{2}(#1)}}
\newcommand{\pg}[1]{\ensuremath{\mathrm{PG}(#1)}}
\newcommand{\gf}[1]{\ensuremath{\mathrm{GF}(#1)}}
\newcommand{\cml}[1]{\ensuremath{\mathit{CM}_{\hspace{-1.5pt}#1}}}
\newcommand{\qml}[1]{\ensuremath{\mathit{QM}_{\hspace{-1.5pt}#1}}}
\newcommand{\imod}[1]
{\mbox{\ensuremath{\hspace{3pt}(\operatorname{mod}\hspace{3pt}#1)}}}

\newcommand{\dash}{\nobreakdash-\hspace{0pt}}
\newcommand{\ba}{\backslash}
\newcommand{\cl}{\operatorname{cl}}
\newcommand{\si}{\operatorname{si}}
\newcommand{\mob}{M\"{o}bius}
\newcommand{\mkt}{\ensuremath{M(K_{3,3})}}
\newcommand{\mkf}{\ensuremath{M(K_{5})}}

\newcommand{\ifc}{internally $4$\nobreakdash-\hspace{0pt}connected}

\begin{document}

\title[Excluding Kuratowski graphs]
{Excluding Kuratowski graphs and their duals from
binary matroids}

\author{Dillon Mayhew, Gordon Royle, Geoff Whittle}

\date{}

\begin{abstract}
We consider some applications of our characterisation
of the \ifc\ binary matroids with no \mkt\dash minor.
We characterise the \ifc\ binary matroids with no minor in
$\mcal{M}$, where \mcal{M} is a subset of
$\{\mkt,M^{*}(K_{3,3}),\mkf,M^{*}(K_{5})\}$ that
contains either \mkt\ or $M^{*}(K_{3,3})$.
We also describe a practical algorithm for testing
whether a binary matroid has a minor in \mcal{M}.
In addition we characterise the growth-rate of
binary matroids with no \mkt\dash minor, and we show that a
binary matroid with no \mkt\dash minor has critical exponent over
\gf{2} at most equal to four.
\end{abstract}

\maketitle

\section{Introduction}
\label{intro}

Earlier, we proved the following theorem.

\begin{theorem}[{\cite[Theorem~1.1]{MRW10}}]
\label{thm1}
An \ifc\ binary matroid $M$ has no minor
isomorphic to \mkt\ if and only if $M$ is either:
\begin{enumerate}[label=\textup{(\roman*)}]
\item cographic,
\item isomorphic to a triangular or triadic
\mob\ matroid, or
\item isomorphic to one of $18$ sporadic
matroids.
\end{enumerate}
\end{theorem}

The $18$ sporadic matroids appearing in \Cref{thm1}
have ground sets of cardinality at most $21$, and have rank
at most $11$.
Their matrix representations appear in
Appendix~B of~\cite{MRW10}.
\mob\ matroids are single-element extensions of the cographic
matroids corresponding to two families of graphs:
The \emph{cubic \mob\ ladder} \cml{2n} is
obtained from an even cycle with vertex sequence $v_{0},\ldots, v_{2n-1}$ by
joining each vertex $v_{i}$ to the antipodal vertex $v_{i+n}$.
(Indices are read modulo~$2n$.)
The \emph{quartic \mob\ ladder} \qml{2n+1} is obtained from an
odd cycle with vertex sequence $v_{0},\ldots, v_{2n}$
by joining each vertex $v_{i}$ to the two antipodal vertices
$v_{i+n}$ and $v_{i+n+1}$.
(Indices are read modulo~$2n+1$.)

Let $r>2$ be an integer and let
$\{e_{1},\ldots, e_{r}\}$ be the standard basis of $\gf{2}^{r}$.
For $1 \leq i \leq r-1$ let $a_{i}$ be the sum of
$e_{i}$ and $e_{r}$, and for $1 \leq i \leq r-2$
let $b_{i}$ be the sum of $e_{i}$ and $e_{i+1}$.
Let $b_{r-1}$ be the sum of $e_{1}$, $e_{r-1}$,
and $e_{r}$.
The rank\dash $r$ \emph{triangular \mob\ matroid},
denoted by $\Delta_{r}$, is represented by the set
$\{e_{1},\ldots, e_{r}, a_{1},\ldots, a_{r-1}, b_{1},\ldots, b_{r-1}\}$.
(We also take this set to be the ground set of $\Delta_{r}$.)
Deleting $e_{r}$ from $\Delta_{r}$ produces a copy of
$M^{*}(\cml{2r-2})$.
It is easy to see that if $r \geq 4$, then
$\Delta_{r}$ has $\Delta_{r-1}$ as a minor.

Now let $r \geq 4$ be an even integer, and again let
$\{e_{1},\ldots, e_{r}\}$ be the standard basis of
$\gf{2}^{r}$.
For $1 \leq i \leq r - 2$ let $c_{i}$ be the sum of
$e_{i}$, $e_{i+1}$, and $e_{r}$.
Let $c_{r-1}$ be the sum of $e_{1}$, $e_{r-1}$, and
$e_{r}$.
The rank\dash $r$ \emph{triadic \mob\ matroid}, denoted
by $\Upsilon_{r}$, is represented by the set
$\{e_{1},\ldots, e_{r}, c_{1},\ldots, c_{r-1}\}$.
If $r \geq 4$ is an even integer then
$\Upsilon_{r} \ba e_{r}$ is isomorphic to $M^{*}(\qml{r-1})$.
If $r>4$, then $\Upsilon_{r}$ has $\Upsilon_{r-2}$ as a minor.

This sequel explores various applications of \Cref{thm1}.
If \mcal{M} is a set of binary matroids, then
\ex{\mcal{M}} is the set of binary matroids that
have no minor isomorphic to a member of \mcal{M}.
Throughout this introduction, we let \mcal{M} be some
subset of $\{\mkt,M^{*}(K_{3,3}),\mkf,M^{*}(K_{5})\}$
that contains either \mkt\ or $M^{*}(K_{3,3})$.
First of all, by using \Cref{thm1}, and the classical
results by Hall and Wagner on graphs with no $K_{3,3}$\dash\ or
$K_{5}$\dash minor, we can obtain additional
characterisations.
In \Cref{classes} we list descriptions of the \ifc\ matroids
in \ex{\mcal{M}}.
Thus we characterise the \ifc\ members in
twelve different families of binary matroids.
Only the smallest of these classes has been characterised
before~\cite{QZ04}.

The Graph Minors Project of Robertson and Seymour showed that
there is a polynomial-time algorithm for testing whether a graph
contains a fixed minor \cite{RS95}.
Similarly, the Matroid Minors Project of Geelen, Gerards, and Whittle \cite{GGW07}
is expected to show that the following problem has a
polynomial-time solution for each \gf{q}\dash representable matroid $N$:
Given a matrix, $A$, over the field \gf{q}, decide whether the
matroid $M[A]$ has an $N$\dash minor.
However, the existence proofs of these algorithms
are very non-constructive.
In \Cref{algorithms} we present algorithms that
could actually be implemented.
In particular, we present an algorithm that will
decide whether $M[A]$ has a minor in
\mcal{M}, where $A$ is a matrix over \gf{2}.
The algorithm runs in $O(n^{13})$ steps, where $n$
is the number of columns in $A$.

A very well-known example due to Seymour \cite{Sey81a}
shows that an oracle algorithm for testing whether
a matroid is binary cannot run in polynomial time relative
to the size of the ground set.
As we discuss in \Cref{algorithms}, the same
example shows that there is no polynomial-time
oracle algorithm for testing whether a matroid is
binary with no minor in \mcal{M}.
However, this difficulty vanishes when we
restrict ourselves to \ifc\ matroids:
There is a polynomial-time oracle algorithm that tests
whether an \ifc\ matroid (not necessarily binary)
belongs to \ex{\mcal{M}}.
We conjecture that this is a general phenomenon:
the problems created by Seymour's examples
can be eliminated with higher connectivity.

\begin{conjecture}
\label{recognise-binary}
There is a polynomial-time oracle algorithm for deciding
if an \ifc\ matroid is binary.
\end{conjecture}

This is an ambitious conjecture.
The next is somewhat more modest.

\begin{conjecture}
\label{recognise-subclass}
There is a polynomial-time oracle algorithm for deciding
whether an \ifc\ matroid belongs to any given proper
minor-closed class of binary matroids.
\end{conjecture}

The matroids in Seymour's example are known as binary spikes.
Spikes are a notorious source of difficulty in matroid theory.
A \emph{spike-like flower} of order $n$ in a 
$3$\dash connected matroid $M$ is a partition
$(P_{1},\ldots, P_{n})$ of the ground set of $M$
such that, for every proper subset $J$ of $\{1,\ldots, n\}$ the
partition
$(\cup_{j \in J} P_{j},E(M) - \cup_{j \in J} P_{j})$
is an exact $3$\dash separation of $M$; and,
for all distinct $i$ and $j$ in $\{1,\ldots, n\}$ we have
$r(P_{i} \cup P_{j}) = r(P_{i}) + r(P_{j}) - 1$.
A rank\dash $n$ spike contains a spike-like flower of order $n$.
We believe the existence of large spike-like flowers is at the
heart of the difficulty of recognising binary matroids.
This belief is encapsulated by the next
conjecture, which is a strengthening of
\Cref{recognise-binary}.

\begin{conjecture}
\label{spike-out}
Let $k$ be a fixed positive integer.
There is a polynomial-time oracle algorithm for deciding
if a $3$\dash connected matroid with no spike-like flower
of order $k$ is a binary matroid.
\end{conjecture}

In the final two sections of the paper, we consider growth-rates and
critical exponents.
In \Cref{extremal} we use \Cref{thm1} to
determine the maximum size of a simple rank\dash $r$ binary matroid
with no \mkt\dash minor.
Moreover, we characterise the matroids that obtain this
upper bound.
This completely resolves a question studied by Kung~\cite{Kun87}.
He showed that a simple rank\dash $r$ binary matroid $M$ without
an \mkt\dash minor has at most $10r$ elements.
\Cref{thm15} shows that, in fact,
$|E(M)| \leq 14r/3 - \alpha(r)$,
where $\alpha(r)$ assumes one of three values
depending on the residue of $r$ modulo~$3$.
Any matroid meeting this bound can be obtained by starting with
either \pg{1,2}, \pg{2,2}, or \pg{3,2},
and then repeatedly adding copies of \pg{3,2}
via parallel connections along points.

If $M$ is a matroid, then its characteristic polynomial,
$\chi(M;t)$, is a polynomial in the variable $t$, and naturally
generalises the chromatic polynomial of a graph.
If $M$ is loopless and representable over \gf{q}, then the
critical exponent of $M$ over $q$, denoted
$c(M;q)$, is the smallest positive integer $k$ such that
$\chi(M; q^{k}) \ne 0$.
The material in \Cref{extremal} shows that $|E(M)|\leq 5r(M)$,
for every simple binary matroid, $M$,
with no \mkt\dash minor.
It therefore follows from Lemma 7.5 in \cite{Kun96} that
the critical exponent of such a matroid is at most $5$.
Kung had already shown that the critical exponent
is at most 10 \cite{Kun87}.
In \Cref{critical} we improve these bounds by showing
that any loopless binary matroid with no \mkt\dash minor
has a critical exponent over \gf{2} of at most $4$.
This result cannot be improved: we also characterise
the matroids with critical
exponent exactly equal to $4$:
They are precisely those with a $3$\dash connected
component isomorphic to \pg{3,2}.

\section{Additional classes}
\label{classes}

Kung initiated the study of binary matroids that have no minor
isomorphic to one of the graphic matroids of the
Kuratowski graphs, $K_{3,3}$ and $K_{5}$ \cite{Kun87}.
We extend his programme here.
Recall that if \mcal{M} is a set of binary matroids, then \ex{\mcal{M}}
is the class of binary matroids that have no minors isomorphic
to members of \mcal{M}.
Thus \Cref{thm1} gives a structural characterisation of
\ex{\{\mkt\}}.
Let \mcal{M} be a subset of the collection
$\{\mkt,M^{*}(K_{3,3}),\mkf,M^{*}(K_{5})\}$
such that \mcal{M} contains either \mkt\ or its dual.
There are exactly twelve such subsets, leading
to twelve classes of the form \ex{\mcal{M}}.
By using \Cref{thm1}, and classical results by
Hall~\cite{Hal43} and Wagner~\cite{Wag37},
we obtain characterisations of the \ifc\ matroids in each
of these twelve classes.
To do so, we occasionally need to check whether certain binary matroids
have particular minors.
We accomplish this task by using the matroid capabilities of
the Sage mathematics package
(\url{www.sagemath.org}).
We start with some preliminary lemmas.

\begin{lemma}
\label{lem7}
The triangular and triadic \mob\ matroids have no
\mkf\dash minors.
\end{lemma}

\begin{proof}
Lemma~3.8 of~\cite{MRW10} states that the only
\ifc\ non-cographic minors of \mob\ matroids are
themselves \mob\ matroids.
Thus if a \mob\ matroid had an \mkf\dash minor it
would imply that \mkf\ is a \mob\ matroid.
It is easily seen that this is not the case.
\end{proof}

The next lemma follows from Wagner's characterisation of
graphs with no $K_{5}$\dash minor (see~\cite[Theorem~1.6]{KM07}).

\begin{lemma}
\label{lem1}
If $M$ is an \ifc\ cographic matroid with no
minor isomorphic to $M^{*}(K_{5})$ then either $M = M^{*}(G)$,
where $G$ is a planar graph, or $M$ is isomorphic to
either $M^{*}(K_{3,3})$ or $M^{*}(\cml{8})$.
\end{lemma}

Now to the characterisations of the twelve families.
All sporadic matroids are described in Appendix~B of~\cite{MRW10}.
The next result follows from \Cref{thm1}, \Cref{lem7},
and a simple computer check.

\begin{theorem}
\label{thm2}
An \ifc\ matroid $M$ belongs to
\ex{\{\mkt,\mkf\}}
if and only if $M$ is either:
\begin{enumerate}[label=\textup{(\roman*)}]
\item cographic,
\item isomorphic to a triangular or triadic
\mob\ matroid, or
\item isomorphic to one of the sporadic matroids
$C_{11}$, $C_{12}$, $M_{5,12}^{a}$, $M_{6,13}$,
$M_{7,15}$, $M_{9,18}$, or $M_{11,21}$.
\end{enumerate}
\end{theorem}

It is easy to check that
$\Delta_{6}$ has an $M^{*}(K_{5})$\dash minor, and therefore
$\Delta_{r}$ has an $M^{*}(K_{5})$\dash minor for all $r \geq 6$.
On the other hand $\Delta_{r}$ has no
$M^{*}(K_{5})$\dash minor if $r \in \{3, 4,5\}$.
Similarly $\Upsilon_{r}$ has an $M^{*}(K_{5})$\dash minor
if $r \geq 6$, but $\Upsilon_{4}$ has no
$M^{*}(K_{5})$\dash minor.
We note that $\Delta_{3} \cong F_{7}$ and
$\Upsilon_{4} \cong F_{7}^{*}$.
The next theorem follows from these facts,
and by applying \Cref{thm1}, \Cref{lem1}, and
some computer tests.
The sporadic matroid $T_{12}$ was introduced in
\cite{Kin97}.

\begin{theorem}
\label{thm5}
An \ifc\ matroid $M$
belongs to \ex{\{\mkt, M^{*}(K_{5})\}}
if and only if $M$ is either:
\begin{enumerate}[label=\textup{(\roman*)}]
\item planar graphic,
\item isomorphic to one of the cographic matroids
$M^{*}(K_{3,3})$ or $M^{*}(\cml{8})$,
\item isomorphic to one of the \mob\ matroids
$F_{7}$, $F_{7}^{*}$, $\Delta_{4}$, $\Delta_{5}$, or
\item isomorphic to one of the $18$ sporadic matroids
of \textup{\Cref{thm1}}, other than $T_{12}$.
\end{enumerate}
\end{theorem}

Next we consider \ex{\{\mkt,M^{*}(K_{3,3})\}}.
A result due to Hall~\cite{Hal43} implies that the
only $3$\dash connected cographic matroids with no
$M^{*}(K_{3,3})$\dash minor are $M^{*}(K_{5})$, and cycle matroids of
planar graphs.
The only \mob\ matroids with no $M^{*}(K_{3,3})$\dash minor
are $\Delta_{3}$, $\Upsilon_{4}$, and
$\Upsilon_{6}$.
Next we consider the sporadic matroids.
As $\Delta_{4}$ has an $M^{*}(K_{3,3})$\dash minor, we
need only consider sporadic matroids with no
$\Delta_{4}$\dash minor.
By a result in \cite{Kin97} the matroid $T_{12}$ has a transitive
automorphism group, so $T_{12}\ba e$ and $T_{12}/e$ are well-defined.
Corollary~2.15 of~\cite{MRW10} says that the only
\ifc\ non-cographic matroids in \ex{\{\mkt,\Delta_{4}\}}
are $F_{7}$, $F_{7}^{*}$, \mkf, $T_{12} \ba e$, $T_{12} / e$, and $T_{12}$.
None of these matroids has an $M^{*}(K_{3,3})$\dash minor.
Both $T_{12}$ and $T_{12}/e$ are among the sporadic matroids
of \Cref{thm1}, while $F_{7}\cong \Delta_{3}$, $F_{7}^{*}\cong \Upsilon_{4}$, and
$T_{12}\ba e\cong \Upsilon_{6}$ are all \mob\ matroids.
The next result follows.

\begin{theorem}
\label{thm7}
An \ifc\ matroid $M$
belongs to \ex{\{\mkt, M^{*}(K_{3,3})\}}
if and only if $M$ is either:
\begin{enumerate}[label=\textup{(\roman*)}]
\item planar graphic,
\item isomorphic to the cographic matroid $M^{*}(K_{5})$,
\item isomorphic to one of the \mob\ matroids
$F_{7}$, $F_{7}^{*}$, or $T_{12}\ba e$, or
\item isomorphic to one of the sporadic matroids
\mkf, $T_{12} / e$, or $T_{12}$.
\end{enumerate}
\end{theorem}

The next theorems are easy consequences of results
stated above.

\begin{theorem}
\label{thm8}
An \ifc\ matroid $M$ belongs to
\ex{\{\mkt, \mkf, M^{*}(K_{5})\}} if and only if
$M$ is either:
\begin{enumerate}[label=\textup{(\roman*)}]
\item planar graphic,
\item isomorphic to one of the cographic matroids
$M^{*}(K_{3,3})$ or $M^{*}(\cml{8})$,
\item isomorphic to one of the \mob\ matroids
$F_{7}$, $F_{7}^{*}$, $\Delta_{4}$, $\Delta_{5}$, or
\item isomorphic to one of the sporadic matroids
$C_{11}$, $C_{12}$, $M_{5,12}^{a}$, $M_{6,13}$,
$M_{7,15}$, $M_{9,18}$, or $M_{11,21}$.
\end{enumerate}
\end{theorem}

\begin{theorem}
\label{thm10}
An \ifc\ matroid $M$ belongs to
\ex{\{\mkt, M^{*}(K_{3,3}), \mkf\}} if and only if
$M$ is either:
\begin{enumerate}[label=\textup{(\roman*)}]
\item planar graphic,
\item isomorphic to the cographic matroid $M^{*}(K_{5})$, or
\item isomorphic to one of the \mob\ matroids
$F_{7}$, $F_{7}^{*}$, or $T_{12}\ba e$.
\end{enumerate}
\end{theorem}

Finally, we have the following characterisation, which has
already been proved by Qin and Zhou~\cite{QZ04}.

\begin{theorem}
\label{thm12}
An \ifc\ matroid $M$ belongs to
\ex{\{\mkt, M^{*}(K_{3,3}), \mkf, M^{*}(K_{5})\}} if and
only if $M$ is either:
\begin{enumerate}[label=\textup{(\roman*)}]
\item planar graphic, or
\item isomorphic to one of $F_{7}$ or $F_{7}^{*}$.
\end{enumerate}
\end{theorem}

\Cref{thm1} and \Crefrange{thm2}{thm12} gives us
characterisations of seven families.
Dualising these theorems gives us five additional
characterisations
(since \Cref{thm7,thm12} characterise self-dual classes).

\section{Polynomial-time algorithms}
\label{algorithms}

Let \mcal{M} be a set of binary matroids.
We consider the following computational problem.
\begin{quote}
\textsc{Membership of \ex{\mcal{M}}}\\
\textsc{Input:} A \gf{2} matrix representing the matroid, $M$.\\
\textsc{Question:} Does $M$ belong to the class \ex{\mcal{M}}?
\end{quote}
The fact that this problem has a polynomial-time solution is expected
to follow from the Matroid Minors Project of Geelen, Gerards, and Whittle
(see \cite{GGW07}).
However, the proofs in that project are highly non-constructive,
and the algorithms that follow from them are not implementable.
In this section we will describe a practical algorithm for
solving \textsc{Membership of \ex{\mcal{M}}} when
\mcal{M} is a subset of
$\{\mkt,M^{*}(K_{3,3}),\mkf,M^{*}(K_{5})\}$
that contains \mkt\ or $M^{*}(K_{3,3})$.
We start with some preliminary material.

The symmetric difference of sets $Z_{1}$ and $Z_{2}$ is
denoted by $Z_{1} \triangle Z_{2}$.
Suppose that $M$ is a binary matroid.
A \emph{cycle} of $M$ is either the empty set,
or a set that can be partitioned into circuits.
Binary matroids are characterised by the fact that the
symmetric difference of any two cycles is another cycle
\cite[Theorem~9.1.2]{Oxl11}.

Let $M_{1}$ and $M_{2}$ be two binary matroids on the
ground sets $E_{1}$ and $E_{2}$ respectively.
Let \mcal{Z} be the collection
\[
\{Z_{1} \triangle Z_{2} \colon
Z_{i}\ \text{is a cycle of}\ M_{i}\ \text{for}\ i=1, 2,\
\text{and}\ Z_{1}\cap E_{2}=Z_{2}\cap E_{1}\}.
\]
Then \mcal{Z} is the collection of cycles of a binary
matroid on the ground set $E_{1} \triangle E_{2}$
(see \cite{Sey80} or \cite[Lemma~9.3.1]{Oxl11}).
We denote this matroid $M_{1} \triangle M_{2}$.

\begin{proposition}[{\cite[(4.4)]{Sey80}}]
\label{prop1}
Suppose that $M_{1}$ and $M_{2}$ are binary matroids on
the sets $E_{1}$ and $E_{2}$ respectively.
If $I$ and $J$ are disjoint subsets of $E_{1} - E_{2}$
then $(M_{1} \triangle M_{2}) / I \ba J = (M_{1} / I \ba J) \triangle M_{2}$.
\end{proposition}

Next we define $1$\dash, $2$\dash, and $3$\dash sums
of binary matroids, following the path taken by
Seymour \cite{Sey80}.
This means that the sums defined here are
slightly different from those in \cite{Oxl11}.
If $E_{1}$ and $E_{2}$ are disjoint, and neither
$E_{1}$ nor $E_{2}$ is empty, then
$M_{1} \triangle M_{2}$ is the \emph{$1$\dash sum} of
$M_{1}$ and $M_{2}$, denoted $M_{1} \oplus_{1}M_{2}$.
If $E_{1}\cap E_{2}=\{p\}$, where
$p$ is neither a loop nor a coloop in $M_{1}$ or $M_{2}$, and
$|E_{1}|, |E_{2}| \geq 3$, then $M_{1} \triangle M_{2}$ is the
\emph{$2$\dash sum} of $M_{1}$ and $M_{2}$, denoted
$M_{1} \oplus_{2} M_{2}$.
We say that $p$ is the \emph{basepoint} of the $2$\dash sum.
Finally, suppose that $E_{1} \cap E_{2} = T$ and
assume that the following conditions hold:
\begin{enumerate}[label=\textup{(\roman*)}]
\item $T$ is a triangle in both $M_{1}$ and $M_{2}$,
\item $T$ contains a cocircuit in neither $M_{1}$ nor $M_{2}$, and
\item $|E_{1}|, |E_{2}| \geq 7$.
\end{enumerate}
In this case $M_{1} \triangle M_{2}$ is the
\emph{$3$\dash sum} of $M_{1}$ and $M_{2}$, denoted
$M_{1} \oplus_{3} M_{2}$.

\begin{proposition}[{\cite[(2.1)]{Sey80}}]
\label{county}
If $(X_{1},X_{2})$ is a $1$\dash separation of the binary matroid $M$,
then $M=(M|X_{1})\oplus_{1}(M|X_{2})$.
Conversely, if $M=M_{1}\oplus_{1}M_{2}$, then
$(E(M_{1}),E(M_{2}))$ is a $1$\dash separation of $M$.
\end{proposition}

The next result is easy, and also follows from \cite[Proposition~4.2.20]{Oxl11}
and \cite[(2.1)]{Oxl11}.

\begin{proposition}
\label{fender}
Let $M=M_{1}\oplus_{1}M_{2}$ and let $N$ be a connected
matroid.
Then $M$ has an $N$\dash minor if and only if
either $M_{1}$ or $M_{2}$ has an $N$\dash minor.
\end{proposition}

\begin{proposition}[{\cite[(2.6)]{Sey80}}]
\label{prop4}
If $(X_{1},X_{2})$ is an exact $2$\dash separation
of the binary matroid $M$,
then there are binary matroids $M_{1}$ and $M_{2}$
on the ground sets $X_{1} \cup p$ and $X_{2} \cup p$,
where $p \notin X_{1} \cup X_{2}$, such that
$M = M_{1} \oplus_{2} M_{2}$.
Conversely, if $M = M_{1} \oplus_{2} M_{2}$
then $(E(M_{1}) - E(M_{2}), E(M_{2})  - E(M_{1}))$
is a $2$\dash separation of $M$.
\end{proposition}

We can deduce the next result from
\cite[Proposition~8.3.5]{Oxl11} and \cite[(2.6)]{Sey80}.

\begin{proposition}
\label{falcon}
Let $M=M_{1}\oplus_{2} M_{2}$ and let $N$ be a
$3$\dash connected matroid.
Then $M$ has an $N$\dash minor if and only if
either $M_{1}$ or $M_{2}$ has an $N$\dash minor.
\end{proposition}

\begin{proposition}[{\cite[(2.9)]{Sey80}}]
\label{prop5}
Suppose that $(X_{1}, X_{2})$ is an exact $3$\dash separation
of the binary matroid $M$ such that $\min\{|X_{1}|, |X_{2}|\} \geq 4$.
Then there are binary matroids $M_{1}$ and $M_{2}$ on the ground
sets $X_{1} \cup T$ and $X_{2} \cup T$ respectively, where $T$ is
disjoint from $X_{1} \cup X_{2}$, such that
$M = M_{1} \oplus_{3} M_{2}$. 
Conversely, if $M = M_{1} \oplus_{3} M_{2}$,
then $(E(M_{1}) - E(M_{2}), E(M_{2}) - E(M_{1}))$ is an
exact $3$\dash separation of $M$.
\end{proposition}

\begin{proposition}[{\cite[(4.1)]{Sey80}}]
\label{prop10}
Let $M$ be the binary matroid $M_{1}\oplus_{3} M_{2}$.
If $M$ is $3$\dash connected then $M_{1}$ and $M_{2}$ are
isomorphic to minors of $M$.
\end{proposition}

\begin{proposition}
\label{lem3}
Let $M$ be the binary matroid $M_{1}\oplus_{3} M_{2}$.
Assume $M$ is $3$\dash connected and let $N$ be
an \ifc\ binary matroid such that $|E(N)|\geq 4$ and
$N$ has no triad.
Then $M$ has an $N$\dash minor if and only if either
$M_{1}$ or $M_{2}$ has an $N$\dash minor.
\end{proposition}

\begin{proof}
Let $E_{1}$ and $E_{2}$ be the ground sets of $M_{1}$ and $M_{2}$
respectively, so that $|E_{1}|,|E_{2}|\geq 7$.
We will assume that $E_{1}\cap E_{2}=T$, where
$T$ is a coindependent triangle in both $M_{1}$ and $M_{2}$.
The `if' direction of the proof follows from \Cref{prop10}.
To prove the `only if' direction, we assume that neither
$M_{1}$ nor $M_{2}$ has an $N$\dash minor, and yet $M$ does.
 Amongst such counterexamples, assume $M$ has been chosen so
 that $|E(M)|$ is as small as possible.

It cannot be the case that $M$ is isomorphic to $N$, or else
\Cref{prop5} would imply that $N$ is not \ifc.
Therefore $M$ has a proper $N$\dash minor.
Furthermore, $N$ is not a wheel, since
all wheels have triads.
Therefore we can apply Seymour's Splitter Theorem
\cite[(7.3)]{Sey80}.
There is a $3$\dash connected minor, $M'$, of $M$
such that $M'$ has an $N$\dash minor, and
$|E(M)|-|E(M')|=1$.
Let $e$ be the element in $E(M)-E(M')$.
Without loss of generality, we can assume that
$e$ is in $E_{1}-T$.
If $M'=M\ba e$, then let $M_{1}'=M_{1}\ba e$, and
if $M'=M/e$, let $M_{1}'=M_{1}/e$.
\Cref{prop1} implies $M'=M_{1}'\triangle M_{2}$.
Since neither $M_{1}'$ nor $M_{2}$ has an $N$\dash minor,
and yet $M'$ does, it follows that
$M_{1}'\triangle M_{2}$ is not the $3$\dash sum of
$M_{1}'$ and $M_{2}$, or else the minimality of
$M$ is contradicted.
Therefore, either $T$ is not a triangle in
$M_{1}'$, or $T$ contains a cocircuit in $M_{1}'$, or
$|E(M_{1}')|<7$.
We eliminate these possibilities one by one.

\begin{sublemma}
\label{sizzle}
$T$ is a triangle in $M_{1}'$.
\end{sublemma}

\begin{proof}
If $T$ is not a triangle, then $M_{1}'$ must be $M_{1}/e$,
and $e$ must be parallel to an element, $x\in T$, in $M_{1}$.
Let $C$ be a circuit of $M_{2}$ such that $C\cap T=\{x\}$
($C$ exists because $T$ is coindependent in $M_{2}$).
It is easy to see that $(C-x)\cup e$ is a circuit of $M=M_{1}\triangle M_{2}$.
Now $(E_{1}-T,E_{2}-T)$ is a $3$\dash separation of
$M$, by \Cref{prop5}, and $e$ is in $E_{1}\cap \cl_{M}(E_{2})$.
Thus $(E_{1}-(T\cup e),E_{2}-T)$ is a $2$\dash separation of
$M'=M/e$, and this contradicts the fact that $M'$ is $3$\dash connected.
\end{proof}

\begin{sublemma}
\label{jockey}
$T$ does not contain a cocircuit in $M_{1}'$.
\end{sublemma}

\begin{proof}
Certainly $T$ does not contain a cocircuit in $M_{1}$.
If it contains a cocircuit in $M_{1}'$, then $M_{1}'=M_{1}\ba e$,
and there is a cocircuit, $C_{1}^{*}$, of $M_{1}$ such that
$C_{1}^{*}\subseteq T\cup e$ and $e\in C_{1}^{*}$.
The intersection $T\cap C_{1}^{*}$ contains exactly
two elements \cite[Theorem~9.1.2]{Oxl11}.
Let $T=\{x_{1},x_{2},x_{3}\}$, and assume that
$T\cap C_{1}^{*}=\{x_{1},x_{2}\}$.

In $M_{2}$, consider a basis, $B$, that contains $x_{2}$ and $x_{3}$.
Then $\cl_{M_{2}}(B-x_{2})$ is a hyperplane that intersects
$T$ exactly in $x_{3}$.
Thus there is a cocircuit, $C_{2}^{*}$, of $M_{2}$ such that
$C_{2}^{*}\cap T=\{x_{1},x_{2}\}$.
Now $M_{1}^{*}\triangle M_{2}^{*}=M^{*}$
(see \cite[p.~319]{Sey80}).
From this we can deduce that $(C_{2}^{*}-\{x_{1},x_{2}\})\cup e$ is a
cocircuit in $M$.
Hence $(E_{1}-T,E_{2}-T)$ is a $3$\dash separation of
$M$, and $e$ is in $E_{1}\cap \cl_{M}^{*}(E_{2})$.
Therefore $(E_{1}-(T\cup e),E_{2}-T)$ is a $2$\dash separation
of $M'=M\ba e$, a contradiction.
\end{proof}

By~\ref{sizzle} and~\ref{jockey}, we must now assume that $|E(M_{1}')|<7$, and
hence $|E_{1}|=7$.
From \cite[(4.3)]{Sey80}, we know that $M_{1}$ is $3$\dash connected,
except that there may exist parallel classes of size two that contain
elements of $T$.
In particular, $M_{1}$ contains no series pair and no coloop.
Any $7$\dash element binary matroid with rank at least $4$
that contains a triangle also contains a series pair or coloop,
so $2\leq r(M_{1}')\leq r(M_{1})\leq 3$.

Assume $r(M_{1}')=3$, so the complement
of $T$ in $M_{1}'$ is a cocircuit of size at
most three.
Let $C^{*}$ be this cocircuit.
From the fact that $(M_{1}')^{*}\triangle M_{2}^{*}=(M')^{*}$,
we can see that $C^{*}$ is a cocircuit of $M'$.
Since $N$ has no triad, there is an
element, $x\in C^{*}$, such that
$M'/x=(M_{1}'/x)\triangle M_{2}$ has an $N$\dash minor.
Therefore either $r(M_{1}')=2$ and
$M_{1}'\triangle M_{2}$ has an $N$\dash minor, or
$r(M_{1}'/x)=2$, and $(M_{1}'/x)\triangle M_{2}$
has an $N$\dash minor.
In either case it is easy to see that
$M_{1}'\triangle M_{2}$ or
$(M_{1}'/x)\triangle M_{2}$
is obtained from $M_{2}$ by possibly
deleting elements of $T$ and adding parallel
elements to elements of $T$.
As $N$ has no parallel pairs, it now follows that
$M_{2}$ has an $N$\dash minor.
This is a contradiction that completes the proof of the
\namecref{lem3}.
\end{proof}
Now we prove the main result of this section.

\begin{theorem}
\label{maiden}
Let \mcal{M} be a subset of
$\{\mkt,M^{*}(K_{3,3}),\mkf,M^{*}(K_{5})\}$
that contains \mkt\ or $M^{*}(K_{3,3})$.
There is an algorithm that solves
\textsc{Membership of \ex{\mcal{M}}}
in time bounded by $O(|E(M)|^{13})$.
\end{theorem}

\begin{proof}
Let $n$ be $|E(M)|$.
A representation of $M^{*}$ can be produced in $O(n^{3})$ steps.
Therefore we can replace $M$ with $M^{*}$ if necessary, so we
lose no generality in assuming that \mcal{M} contains $M^{*}(K_{3,3})$.
Henceforth we also assume that $M$ has no loops or coloops.

We sketch the procedure for finding an exact
$k$\dash separation (when $k$ is $1$ or $2$),
or an exact $3$\dash separation with at least $4$ elements
on each side.
A more complete description is in \cite{BC95}.
The algorithm involves considering all pairs of disjoint
$k$\dash element subsets (when $k$ is $1$ or $2$)
or $4$\dash element subsets (when $k=3$).
This requires looping at most $O(n^{8})$ times.
We attempt to extend each such pair to a $k$\dash separation.
This involves examining each remaining element
of $E(M)$ (looping $O(n)$ times), and calculating the
rank of a submatrix for each such element (which can be
done in $O(n^{3})$ steps).
Thus it takes at most $O(n^{12})$ steps to search for a
separation certifying that $M$ can be decomposed via a
$1$\dash, $2$\dash, or $3$\dash sum.

Every loopless rank\dash $r$ binary matroid
can be considered as a multiset of points in the projective
space $P = \pg{r - 1,2}$.
If $X \subseteq E(M)$, we use $\cl_{P}(X)$
to denote the span of $X$ in $P$.
Suppose that $(X_{1}, X_{2})$ is an exact
$k$\dash separation of $M$ for some $k \in \{1, 2, 3\}$
with the property that if $k = 3$ then $|X_{1}|,|X_{2}| \geq 4$
and $r_{M}(X_{1}), r_{M}(X_{2}) \geq 3$.
Let $Z = \cl_{P}(X_{1}) \cap \cl_{P}(X_{2})$, and for
$i = 1, 2$ let $M_{i}$ be the binary matroid represented
by the multiset $X_{i} \cup Z$.
Then $M \cong M_{1} \oplus_{k} M_{2}$.
It follows easily that by solving a system of equations
(which takes $O(n^{3})$ steps),
we can produce representations of
$M_{1}$ and $M_{2}$.

Imagine a binary tree with nodes labelled by matroids.
The root is labelled $M$.
If a node is labelled $M'$,
we are allowed to label the children of that node with
$M_{1}$ and $M_{2}$ if $M'$ can be expressed
as $M'=M_{1}\oplus_{k} M_{2}$
for some $k\in \{1,2,3\}$.
For every node, we assume the decomposition of $M'$ into
$M_{1}\oplus_{k} M_{2}$ has been chosen so that
$k$ is as small as possible.
Even so, the decomposition need not be unique.
Now $|E(\si(M_{1}))|+|E(\si(M_{2}))|\leq |E(\si(M))|+6$.
It is easy to prove by induction that the
binary tree has at most $\max\{1,|E(\si(M))|-6\}$ leaves, and
therefore at most $O(n)$ internal nodes.
Each internal node corresponds to finding a
$k$\dash separation.
It follows that we can construct such a tree in
time bounded by $O(n)(O(n^{12})+O(n^{3}))=O(n^{13})$.
Note that the matroids labelling leaves are
\ifc\ (except that they may contain parallel pairs).

Let $\mcal{M}_{\Delta}$ be the subset of \mcal{M}
containing matroids that have triangles.
By our earlier assumption, $M^{*}(K_{3,3})$ is in
$\mcal{M}_{\Delta}$.
Note that no matroid in $\mcal{M}_{\Delta}$ contains
a triad.
We have insisted that we decompose a matroid along a
$3$\dash sum only if it is $3$\dash connected.
Therefore we can apply \Cref{fender,falcon} and \Cref{lem8}.
From these results we deduce that
$M$ has no minor in $\mcal{M}_{\Delta}$ if and only the simplification
of each of the matroids labelling a leaf has no such minor.
The basic classes of \ifc\ binary matroids with no minor in
$\mcal{M}_{\Delta}$ are described either by a theorem in \Cref{classes},
or the dual of such a theorem.
Therefore we now check that the simplification of
each leaf matroid belongs to one of these basic classes.
There are $O(n)$ leaves.
Producing a representation of a dual can be done in
$O(n^{3})$ steps.
We will show that we can test whether a matroid is
cographic, isomorphic to a sporadic matroid,
or a \mob\ matroid in at most $O(n^{7})$ steps, so
we can complete this part of the algorithm with another
$O(n)O(n^{3})O(n^{7})$ steps.

Assume that $M'$ is the simplification of a leaf matroid.
Checking that $M'$ is isomorphic to a specific
sporadic matroid can be done in constant time.
There is an algorithm running in at most $O(n^{3})$
calls to an independence oracle
which will produce a graph that represents $M'$,
or certify that no such graph exists \cite{BC95}.
This also allows us to check whether $M'$ is cographic or
planar graphic.
Each call to an independence oracle can be simulated
in $O(n^{3})$ operations on the matrix, so the total time required
to check whether $M'$ is graphic, cographic, or
planar graphic, is $O(n^{6})$.

To check if $M'$ is a triangular \mob\ matroid, we
consider each matroid of the form $M'\ba e$, and
produce a graph $G$ (if possible), such that
$M^{*}(G)=M'\ba e$.
We then check each such graph to see if it
is a cubic \mob\ ladder.
We can do this by finding all $4$\dash cycles
(in time $O(n^{4})$), and checking that the
edges lying in exactly one such cycle form a Hamiltonian cycle.
The remaining edges must then join opposite
vertices in the cycle.
Assuming this is the case, we check that $e$
forms a circuit with the set of edges not in the Hamiltonian cycle.
This entire process can be completed in
$O(n)(O(n^{6})+O(n^{4}))=O(n^{7})$ steps.
To check if $M'$ is a triadic \mob\ matroid, we
go through a similar process, except that
we find the $3$\dash cycles of $G$.
The edges that lie in exactly one $3$\dash cycle
must form a Hamiltonian cycle, and $e$ must be in a
circuit with the set of edges not in this cycle.

We have now shown that it is possible to test in $O(n^{13})$ steps
whether $M$ has a minor isomorphic to a member of $\mcal{M}_{\Delta}$.
If $M$ has such a minor, then we halt the algorithm.
If $\mcal{M}=\mcal{M}_{\Delta}$, then again, we can halt.
Therefore we now assume that $M$ has no minor in $\mcal{M}_{\Delta}$
but that $\mcal{M}-\mcal{M}_{\Delta}$ is non-empty.
We next produce a decomposition tree for $M^{*}$,
and we test the matroids corresponding to leaves of this
tree to see whether they belong to the basic classes of
matroids with no minor in $\{N^{*}\colon N\in\mcal{M}\}$.
Note that the matroids in $\{N^{*}\colon N\in \mcal{M}-\mcal{M}_{\Delta}\}$
have no triads.
We can again use the results from earlier in this
section to deduce that $M^{*}$ has no minor in
$\{N^{*}\colon N\in \mcal{M}-\mcal{M}_{\Delta}\}$
if and only if the leaf matroids belong to the basic classes.
We already have assumed that $M^{*}$ has no minor in
$\{N^{*}\colon N\in \mcal{M}_{\Delta}\}$.
Therefore we can complete the algorithm with another
$O(n^{13})$ steps.
\end{proof}

\noindent\textbf{Oracle algorithms.}
Historically, matroid computation has often been
discussed in terms of oracle algorithms.
In this case our computational model is a deterministic
Turing Machine equipped with an oracle which can, in
unit time, return the rank of a specified subset of the ground set.
Let \mcal{M} be a subset of the family
$\{\mkt,M^{*}(K_{3,3}),\mkf,M^{*}(K_{5})\}$
that contains either $M(K_{3,3})$ or $M^{*}(K_{3,3})$.
We now briefly discuss the difficulty of testing membership
in \ex{\mcal{M}} via oracle algorithms.

A well-known example of Seymour's~\cite{Sey81a} shows that
the number of calls to the oracle needed to test
whether a matroid is binary is exponential in terms of the
size of the ground set.
We will use the same example to show that an exponential
number of calls is required to test whether a matroid belongs to
\ex{\mcal{M}}.
For $r \geq 3$ let $\{e_{1},\ldots, e_{r}\}$ be the standard
basis of the vector space over \gf{2} with dimension $r$.
Let $d$ be the sum of $e_{1},\ldots, e_{r}$, and for
$1 \leq i \leq r$ let $d_{i}$ be the sum of $d$ and $e_{i}$.
The binary matroid represented by the set
$\{e_{1},\ldots, e_{r}, d_{1},\ldots, d_{r}\}$
is known as the rank\dash $r$ \emph{binary spike}.
We will denote this matroid by $N_{r}$.
If $H$ is a subset of $E(N_{r})$ such that
$|H \cap \{d_{1},\ldots, d_{r}\}|$ is odd and
$|H \cap \{e_{i}, d_{i}\}| = 1$ for $1 \leq i \leq r$,
then $H$ is a circuit-hyperplane of $N_{r}$.
Let $N_{r}(H)$ be the matroid obtained from $N_{r}$
by relaxing $H$.
It is not difficult to prove by induction on $r$ that
$N_{r}$ has no minor in \mcal{M}.
However, $N_{r}(H)$ is non-binary, so certainly does
not belong to \ex{\mcal{M}}.
In the worst case, an oracle algorithm will have to check
each of the $2^{r-1}$ candidate sets, $H$, to decide whether the
matroid it is considering is isomorphic to $N_{r}$ or $N_{r}(H)$.
Therefore testing whether a matroid belongs to
 \ex{\mcal{M}} requires exponentially many calls to the oracle
 relative to the size of the ground set.

The matroid $N_{r}$ contains
many $3$\dash separations.
If we restrict our attention to \ifc\ matroids the
situation changes dramatically.
Seymour~\cite{Sey81a} shows that there is
an algorithm which, given a matroid $M$ (not necessarily
binary), will either output a graph $G$ such that $M = M(G)$,
or decide that no such graph exists, using a polynomial
number of calls to the oracle.
Using a similar strategy to that in the proof of
\Cref{maiden} we can
decide whether a matroid $M$ is isomorphic to a
\mob\ matroid, using only a polynomial number of calls to an
oracle.
Since it is obviously possible to decide whether a
matroid $M$ is isomorphic to one of a finite number
of sporadic matroids using a constant number of oracle
calls, it follows that we can decide in a polynomial
number of calls to the oracle whether an \ifc\
matroid (not necessarily binary) belongs to \ex{\mcal{M}}.

\section{The growth-rate of $\mathrm{EX}(\{M(K_{3,3})\})$}
\label{extremal}

Kung~\cite{Kun87} investigated
simple rank\dash $r$ matroids of maximum size in \ex{\{\mkt\}}.
He showed that if $N$ and $N'$ are such matroids
and $r(N)=r(N')+1$, then $|E(N)|-|E(N')|\leq 10$.
It then follows by induction that
$|E(M)| \leq 10r(M)$ for any simple matroid $M\in \ex{\{\mkt\}}$.
Using our structure theorem, we show that
$|E(N)|-|E(N')|$ is $4$, $8$, or $2$, depending
on the residue of $r(N)$ modulo $3$.
The average of these three numbers is $14/3$, so we
can prove that $|E(M)|\leq 14r(M)/3$ for any simple matroid $M\in \ex{\{\mkt\}}$.
Moreover, we characterise the simple rank\dash $r$ matroids of maximum
size in the class.

When $r\in\{2,3,4\}$, we define the class
$\mcal{P}_{r}$ to be $\{\pg{r-1,2}\}$.
When $r > 4$ we recursively define $\mcal{P}_{r}$ to be the class
$
\{P(M, \pg{3,2}) \colon M \in \mcal{P}_{r-3}\},
$
where $P(M, \pg{3,2})$ is a parallel connection of
$M$ and $\pg{3,2}$ (see~\cite[Section~7.1]{Oxl11}) along
an arbitrary basepoint.
Note that starting with $r=8$, the class $\mcal{P}_{r}$ contains
non-isomorphic matroids.
It is well known that parallel connections
can be expressed as $2$\dash sums by adding a parallel element.
Therefore the next result follows easily from \Cref{falcon}
and induction.

\begin{proposition}
\label{prop33}
Let $r \geq 2$ be an integer.
If $M$ is in $\mcal{P}_{r}$, then
$M$ has no \mkt\dash minor.
\end{proposition}

For an integer $r\geq 2$, define $h(r)$ to be the
size of matroids in $\mcal{P}_{r}$.
Therefore
\[
h(r)=
\begin{cases}
\frac{14}{3}r-7 &\text{if}\ r\equiv 0\imod{3}\\[0.2ex]
\frac{14}{3}r-\frac{11}{3} &\text{if}\ r\equiv 1\imod{3}\\[0.2ex]
\frac{14}{3}r-\frac{19}{3} &\text{if}\ r\equiv 2\imod{3}
\end{cases}
\]
Let $\alpha(r)$ be $7$, $11/3$, or $19/3$ according to
whether $r$ is equivalent to $0$, $1$, or $2$ modulo $3$.
Thus $h(r)=14r/3-\alpha(r)$.

\begin{theorem}
\label{thm15}
Let $M$ be a simple member of \ex{\{\mkt\}} with rank $r\geq 2$.
Then $|E(M)|\leq h(r)$ and equality holds if and only if
$M\in \mcal{P}_{r}$.
\end{theorem}

\begin{lemma}
\label{lem8}
Assume $M$ is a $3$\dash connected member of \ex{\{\mkt\}} with
rank $r\geq 2$.
Either $|E(M)|\leq 4r-5$, or $M$ is isomorphic to one of
the rank\dash $4$ sporadic matroids
\pg{3,2}, $M_{4,14}$, $M_{4,13}$, $C_{12}$, or $D_{12}$.
\end{lemma}

\begin{proof}
Assume that $M$ is a counterexample with the
smallest possible rank.
Therefore $M$ is a $3$\dash connected
simple matroid in \ex{\{\mkt\}}, with $r\geq 2$ and $|E(M)|>4r-5$, while
$M$ is not isomorphic to any of the five sporadic matroids listed in the statement.
It is easy to see that $r\geq 3$.

Assume that $M$ is \ifc.
All the sporadic matroids from \Cref{thm1} satisfy the
bound $|E(M)|\leq 4r(M)-5$ (with the exceptions of
\pg{3,2}, $M_{4,14}$, $M_{4,13}$, $C_{12}$, and $D_{12}$).
Therefore $M$ is either cographic or isomorphic to a \mob\ matroid.
If $M$ is a \mob\ matroid then $|E(M)|$ is either
$3r - 2$ or $2r - 1$.
If $M$ is cographic then $|E(M)|\leq 3r-3$
\cite[Lemma~14.10.2]{Oxl11}.
Since $r\geq 3$, both $3r-2$ and $2r-1$ are bounded above by $4r-5$.
Thus we have a contradiction in any case, so $M$ is not \ifc.

By \Cref{prop5} we can express $M$ as $M_{1} \oplus_{3} M_{2}$.
Let $T$ be $E(M_{1})\cap E(M_{2})$.
Let $r_{i}$ be the rank of $M_{i}$ for $i=1,2$.
\Cref{prop5} implies $r_{1} + r_{2} - r = 2$.
Certainly $r_{1},r_{2}>2$, or else $T$
contains a cocircuit in $M_{1}$ or $M_{2}$.
Hence $r_{1}, r_{2} < r$.
We see from \cite[(4.3)]{Sey80} that $\si(M_{1})$ and $\si(M_{2})$
are $3$\dash connected.
Moreover $\si(M_{i})$ has no \mkt\dash minor for $i = 1,2$, by \Cref{prop10}.
Therefore the \namecref{lem8} holds for $\si(M_{1})$ and $\si(M_{2})$
by our inductive assumption.

Assume that $\si(M_{2})$ is isomorphic to
\pg{3,2}, $M_{4,14}$, $M_{4,13}$, $C_{12}$, or $D_{12}$.
If $\cl_{M_{1}}(T)$ is not coindependent in $M_{1}$, then $\cl_{M_{1}}(T)$
is $2$\dash separating in $M_{1}$, which 
contradicts \cite[(4.3)]{Sey80}.
Therefore we can deduce that $\si(M_{1})$
has corank at least three.
Now we can apply \cite[Theorem~3.6]{Oxl87a},
which tells us that $M_{1}$ has a minor, $M_{1}'$,
isomorphic to $M(K_{4})$ and containing the triangle $T$.
Now $M_{1}'\triangle M_{2}$ is a minor of $M$, by \Cref{prop1},
and $M_{1}'\triangle M_{2}$ is obtained from $M_{2}$ by performing
a $\Delta\text{-}Y$ operation on the triangle $T$.
It follows that $M$ has, as a minor, a matroid
obtained from
\pg{3,2}, $M_{4,14}$, $M_{4,13}$, $C_{12}$, or $D_{12}$
by performing a $\Delta\text{-}Y$ operation.
But any such matroid contains a \mkt\dash minor,
by results from \cite[Appendix~C]{MRW10}.
From this contradiction and induction we deduce that
$|E(\si(M_{2}))|\leq 4r_{2}-5$.
Symmetrically, $|E(\si(M_{1}))|\leq 4r_{1}-5$.

The only parallel classes of $M_{i}$ have size two
and contain an element of $T$ \cite[(4.3)]{Sey80}.
Note that no element in $T$ can be in a parallel pair
in both $M_{1}$ and $M_{2}$, for that would imply that
$M = M_{1} \oplus_{3} M_{2}$ has a parallel pair.
Let $m$ be the number of elements in $T$ that are
contained in a parallel pair in either $M_{1}$ or $M_{2}$.
Then
$|E(M)|=|E(\si(M_{1}))|+|E(\si(M_{2}))|-m-2(3-m)
\leq |E(\si(M_{1}))|+|E(\si(M_{2}))|-3$.
Thus
\begin{linenomath}
\begin{multline*}
|E(M)| \leq |E(\si(M_{1}))| + |E(\si(M_{2}))|-3\\ 
\leq (4r_{1}-5) + (4r_{2}-5) - 3 = 4(r_{1} + r_{2} - 2) -5 = 4r-5
\end{multline*}
\end{linenomath}
and $M$ is not a counterexample after all.
\end{proof}

\begin{proof}[Proof of \textup{\Cref{thm15}}]
Assume that $M$ is a counterexample with rank $r$, where
$r$ is as small as possible.
Thus $M$ is a simple member of \ex{\{\mkt\}}, and
either $|E(M)|>h(r)$, or $|E(M)|=h(r)$ and $M$ 
does not belong to $\mcal{P}_{r}$.
Certainly $M$ is no larger than
the projective geometry \pg{r-1,2}, so the
result holds if $r\leq 4$.
Hence $r > 4$, and therefore $14r/3-7>4r-5$.
As $|E(M)|\geq h(r)\geq 14r/3-7$,
\Cref{lem8} now implies $M$ is not $3$\dash connected.

Assume that $M = M_{1} \oplus_{1}M_{2}$,
so $M_{1}$ and $M_{2}$ belong \ex{\{\mkt\}}
by \Cref{fender}.
Suppose that $r_{i} = r(M_{i})$ for $i = 1, 2$, so that
$r = r_{1} + r_{2}$.
Since $M$ is simple, $r_{i} > 0$ and hence $r_{i} < r$ for
$i = 1, 2$.
Therefore we can apply the inductive hypothesis and conclude
that
\begin{linenomath}
\begin{multline*}
|E(M)| = |E(M_{1})| + |E(M_{2})| \leq (14r_{1}/3 - \alpha(r_{1})) +
(14r_{2}/3-\alpha(r_{2}))\\
=14r/3-(\alpha(r_{1})+\alpha(r_{2})).
\end{multline*}
\end{linenomath}
But $\alpha(r_{1}) + \alpha(r_{2}) > 7$, regardless of the residue
classes of $r_{1}$ and $r_{2}$ modulo~$3$, so
$|E(M)| < 14r/3 - 7\leq h(r)$, contradicting our earlier statement.
Therefore $M$ is connected, but not $3$\dash connected.

Now $M$ can be expressed as $M_{1} \oplus_{2} M_{2}$.
Let $p$ be the basepoint of the $2$\dash sum.
Let $r_{i} = r(M_{i})$ for $i = 1, 2$.
By \Cref{prop4} we see $r_{1} + r_{2} - r = 1$.
As $M$ has no parallel pairs, $r_{1}, r_{2}> 1$, so $r_{1}, r_{2} < r$.
By \Cref{falcon}, neither $M_{1}$ nor $M_{2}$ has an \mkt\dash minor,
so we can apply the inductive hypothesis to
$\si(M_{1})$ and $\si(M_{2})$.
We can assume that either $M_{1}$ or $M_{2}$ is non-simple,
since otherwise we could add a parallel element to $p$ in $M_{1}$,
and obtain a simple matroid that has one more element than
$M$ despite having no \mkt\dash minor
(since adding parallel elements and taking a $2$\dash sum cannot create
a \mkt\dash minor).
However, it cannot be the case that both $M_{1}$ and $M_{2}$
are non-simple, for then $M$ would be non-simple.
Therefore $|E(M)|=|E(\si(M_{1}))|+|E(\si(M_{2}))|-1$,

Assume that $r_{1},r_{2}\equiv 0\imod{3}$, so
$r\equiv 2\imod{3}$ and
\begin{linenomath}
\begin{multline*}
h(r)=\tfrac{14}{3}r-\tfrac{19}{3}\leq |E(M)|=|E(\si(M_{1}))|+|E(\si(M_{2}))|-1\\
\leq \tfrac{14}{3}r_{1}+\tfrac{14}{3}r_{2}-15
=\tfrac{14}{3}(r_{1}+r_{2}-1)-\tfrac{31}{3}=\tfrac{14}{3}r-\tfrac{31}{3},
\end{multline*}
\end{linenomath}
which is impossible.
We reach a similar contradiction if
the residues of $r_{1}$ and $r_{2}$ are
$(0,2)$, $(2,0)$, or $(2,2)$,
so at least one of $r_{1}$ and $r_{2}$
is equivalent to $1$ modulo $3$.

We consider the case that $r_{1}\equiv 1\imod{3}$ and $r_{2}\equiv 0\imod{3}$,
so that $r\equiv 0\imod{3}$ and
\begin{linenomath}
\begin{multline*}
h(r)=\tfrac{14}{3}r-7\leq |E(M)|=|E(\si(M_{1}))|+|E(\si(M_{2}))|-1\\
\leq \tfrac{14}{3}r_{1}+\tfrac{14}{3}r_{2}-\tfrac{11}{3}-7-1
=\tfrac{14}{3}(r_{1}+r_{2}-1)-\tfrac{21}{3}=\tfrac{14}{3}r-\tfrac{21}{3}.
\end{multline*}
\end{linenomath}
As equality holds throughout, we deduce that
$\si(M_{1})$ and $\si(M_{2})$ belong to
$\mcal{P}_{r_{1}}$ and $\mcal{P}_{r_{2}}$
respectively.
It is easy to see that the parallel connection is an
associative operation.
As $\si(M_{1})$ is formed by taking the parallel connection
of multiple copies of \pg{3,2}, and
$\si(M_{2})$ is similarly formed from copies
of \pg{3,2} and a single copy of \pg{2,2},
it now follows that $M$ belongs to $\mcal{P}_{r}$, so $M$
is not a counterexample after all.
In all the other possible cases, we reach a contradiction in exactly the same way.
\end{proof}

\section{Critical exponents}
\label{critical}

Let $M$ be a loopless rank\dash $r$
\gf{q}\dash representable matroid.
Then $M$ can be considered as a multiset of points
in the projective geometry \pg{r-1,q}.
The \emph{critical exponent} of $M$ over $q$,
denoted by $c(M; q)$, is the smallest
integer $k$ such that there is a
set of hyperplanes, $H_{1},\ldots, H_{k}$, in \pg{r-1,q}
with the property that $H_{1} \cap \cdots \cap H_{k}$ contains
no points of $E(M)$.

The critical exponent depends only on $M$ and $q$,
and not on the particular representation chosen.
We can deduce this fact from a formulation in terms
of the \emph{characteristic polynomial}, denoted $\chi(M; t)$.
Assume $M$ is a matroid on the ground set $E$.
Then
\[
\chi(M; t) = \sum_{A \subseteq E}(-1)^{|A|}t^{r(M) - r(A)}.
\]
Now $\chi(M; q^{k})$ is the number of $k$\dash tuples of
hyperplanes, $(H_{1},\ldots, H_{k})$, in \pg{r-1,q}
satisfying $H_{1}\cap\cdots \cap H_{k}\cap E(M)=\emptyset$
\cite[Theorem~4.1]{Kun96}.
From this it follows that $\chi(M; q^{k})\geq 0$ for all
positive integers $k$, and $c(M; q)$ is the least positive
integer $k$ such that $\chi(M; q^{k})>0$.
Note that, if $k\geq c(M;q)$, then
$\chi(M;q^{k})>0$.
It is obvious that if $M$ has a loop, then $\chi(M;t)$ is
identically zero.
If $e\in E(M)$, then $c(M \ba e; q) \leq c(M; q)$.

Kung~\cite{Kun87} looked at the critical exponent over \gf{2}
of binary matroids with no \mkt\dash minor.
He showed that if $M \in \ex{\{\mkt\}}$ is loopless then
$c(M; 2) \leq 10$.
By using \Cref{thm15} as well as \cite[Lemma~3.1]{Kun86},
we can improve this to $c(M; 2)\leq 5$, since
$|E(M)|\leq 5r(M)$ for every simple matroid $M \in \ex{\{\mkt\}}$.
In this section we improve this further to $c(M; 2) \leq 4$,
and show that this bound cannot be improved.
In particular, we show that if $c(M; 2) = 4$, then $M$ has a
$3$\dash connected component isomorphic to \pg{3,2}.

\begin{lemma}
\label{lem11}
Let $M$ be an \ifc\ binary matroid with no \mkt\dash minor.
Then $c(M; 2) \leq 4$, and if $c(M; 2) = 4$,
then $M$ is isomorphic to \pg{3,2}.
\end{lemma}

\begin{proof}
It is easy to see that \pg{3,2} has critical
exponent~$4$ over \gf{2} (see~\cite[Section~8.1]{Kun96})
so we let $M$ be an \ifc\ member of \ex{\{\mkt\}} other
than \pg{3,2}.
Assume that $M = M^{*}(G)$ for some graph $G$.
Because $M$ is connected, $G$ has no isthmus.
Jaeger showed that $G$ has a nowhere-zero
$8$\dash flow~\cite{Jae76}.
The number of such flows is $\chi(M; 8)$ \cite[Theorem~4.6]{Kun96}.
Hence $\chi(M; 8)> 0$ and thus $c(M;2) \leq 3$.
Consider the \gf{2}\dash representation of
$\Delta_{r}$ discussed in \Cref{intro}.
Each point of \pg{r - 1,2} corresponds to a vector $(x_{1},\ldots, x_{r})$.
Let $H_{1}$, $H_{2}$, and $H_{3}$ be the hyperplanes of
\pg{r - 1,2} defined, respectively, by the
equations $x_{r} = 0$, $x_{1} + \cdots + x_{r} = 0$, and
$\sum x_{i}=0$, where the final sum is taken over all
odd indices in $\{1,\ldots, r-1\}$.
It is easy to see that no point of $M$
is contained in $H_{1}\cap H_{2}\cap H_{3}$,
so $c(\Delta_{r}; 2) \leq 3$.
Similarly,
no point of $\Upsilon_{r}$ is contained in the
hyperplane defined by $x_{1} + \cdots + x_{r} = 0$,
so $c(\Upsilon_{r}; 2) \leq 1$.

Now we can assume that $M$ is neither cographic
nor a \mob\ matroid, so
$M$ is isomorphic to one
of the sporadic matroids in \Cref{thm1}.
The largest such matroid with rank~$4$ is
\pg{3,2}, and it known that every proper minor of
this matroid has critical exponent at most three over
\gf{2}~\cite[Section~8.1]{Kun96}.
Thus we will assume that $r(M) \geq 5$.
The sporadic matroid $T_{12}$ has rank~$6$.
By examining the matrix representation of
$T_{12}$ in~\cite[Appendix~B]{MRW10}, we see that
no point of $T_{12}$ is contained
in the hyperplane defined by $x_{1} + \cdots + x_{6} = 0$.
Thus $c(T_{12}; 2) \leq 1$.
Let $A$ be the matrix in~\cite[Appendix~B]{MRW10}
such that $[I_{5}|A]$ represents the rank\dash $5$ sporadic
matroid $M_{5,12}^{a}$.
If
\[
H_{5,12}^{a} =
\left[
\begin{array}{ccccc}
1 & 1 & 0 & 0 & 0\\
1 & 0 & 1 & 0 & 1\\
1 & 1 & 1 & 1 & 1
\end{array}
\right]
\]
then the matrix product $H_{5,12}^{a}[I_{5}|A]$ contains no
zero columns.
This means that no point of $M_{5,12}^{a}$ is
contained in all three of the hyperplanes defined
by $x_{1}+x_{2} = 0$, $x_{1} + x_{3} + x_{5} = 0$,
and $x_{1} + x_{2} + x_{3} + x_{4} + x_{5} = 0$.
Thus $c(M_{5,12}^{a}; 2) \leq 3$.

In the same way we can show that $M_{5,13}$, $M_{6,13}$,
$M_{7,15}$, $M_{9,18}$, and $M_{11,21}$ all have critical
exponent at most three by examining the matrices
\[
H_{5,13}=
\begin{bmatrix}
1& 0& 0& 0& 1\\
1& 0& 1& 0& 1\\
1& 1& 1& 1& 0
\end{bmatrix},\ 
H_{6,13}=
\begin{bmatrix}
0& 1& 0& 0& 0& 0\\
1& 0& 0& 1& 1& 0\\
1& 1& 1& 1& 1& 1 
\end{bmatrix},
\]
\[
H_{7,15}=
\begin{bmatrix}
1& 0& 1& 0& 1& 0& 1\\
1& 1& 1& 1& 1& 1& 1 
\end{bmatrix},\  
H_{9,18}=
\begin{bmatrix}
1& 1& 0& 0& 0& 0& 0& 0& 0\\
1& 0& 1& 0& 1& 0& 1& 0& 1\\
1& 1& 1& 1& 1& 1& 1& 1& 1
\end{bmatrix},
\]
and
\[
H_{11,21} =
\left[
\begin{array}{ccccccccccc}
0&0&0&0&0&0&1&0&0&1&0\\
1&0&1&0&1&0&0&0&1&1&1\\
1&1&1&1&1&1&1&1&1&1&1
\end{array}
\right].
\]

As every sporadic matroid in \Cref{thm1} can be produced
from one of \pg{3,2}, $M_{5,12}^{a}$, 
$M_{5,13}$, $T_{12}$, $M_{6,13}$, $M_{7,15}$, $M_{9,18}$, or $M_{11,21}$
by deleting elements, the proof is complete.
\end{proof}

Next we come to the main result of this section.

\begin{theorem}
\label{thm19}
Let $M$ be a loopless binary matroid with no \mkt\dash minor.
Then $c(M;2) \leq 4$, and if $c(M; 2) = 4$,
then either
\begin{enumerate}[label=\textup{(\roman*)}]
\item $M$ is isomorphic to \pg{3,2}, or
\item $M$ can be expressed as the $1$\dash\ or $2$\dash sum
of $M_{1}$ and $M_{2}$, where $M_{1}, M_{2}$ belong to
\ex{\{\mkt\}}, and either $c(M_{1}; 2) = 4$ or
$c(M_{2}; 2) = 4$.
\end{enumerate}
\end{theorem}

\begin{proof}
Let $M$ be a minor-minimal counterexample
to the theorem.
\Cref{lem11} shows that $M$ cannot be \ifc.
Assume that $M$ is not connected, so that
$M$ can be expressed as $M_{1} \oplus_{1}M_{2}$.
Clearly $M_{1}$ and $M_{2}$ are loopless members of \ex{\{\mkt\}}.
It is well known, and easy to verify, that
$\chi(M; t) = \chi(M_{1}; t)\chi(M_{2}; t)$.
By the inductive hypothesis, $c(M_{i};2)\leq 4$,
meaning that $\chi(M_{i}; 16) > 0$ for $i=1,2$.
Hence $\chi(M; 16) > 0$, so $c(M;2)\leq 4$.
Since $M$ is a counterexample, $c(M;2)=4$, so
$\chi(M; 8) = 0$.
Therefore $\chi(M_{i};8)=0$ for some $i\in\{1,2\}$.
This implies that $c(M_{i}; 2) = 4$ for some
$i \in \{1, 2\}$.
However, $M$ now satisfies statement~(ii)
so it is not a counterexample at all.

Now we must assume that $M$ is connected.
Assume $M$ can be expressed as
$M_{1}\oplus_{2} M_{2}$, where
$p$ is the basepoint of the $2$\dash sum.
Again, $M_{1}$ and $M_{2}$ are loopless members of \ex{\{\mkt\}}.
Walton and Welsh~\cite[(7)]{WW80} note the following
relation:
\begin{equation}
\label{eqn7}
\chi(M; t) = \frac{\chi(M_{1}; t)\chi(M_{2}; t)}{t - 1}
+ \chi(M_{1} / p; t)\chi(M_{2} / p; t).
\end{equation}
Note that if $M_{i}/p$ is loopless, then, by earlier
discussion, $\chi(M_{i}/p;k)\geq 0$ for all positive integers $k$.
The same statement holds if $M_{i}/p$ has a loop, for then
$\chi(M_{i}/p;t)$ is identically zero.
Since $M_{1}$ and $M_{2}$ are isomorphic to proper minors of
$M$ it follows that $\chi(M_{i}; 16) > 0$ for $i = 1, 2$.
Now~\eqref{eqn7} implies that $\chi(M; 16) > 0$, so
$c(M; 2) \leq 4$.
Therefore it must be the case that $c(M; 2) = 4$, so
that $\chi(M; 8) = 0$.
It follows that either
$\chi(M_{1}; 8) = 0$ or $\chi(M_{2}; 8) = 0$.
Then $M$ satisfies statement~(ii) of the theorem, and
we again have a contradiction.

Finally, we assume that $M$ is $3$\dash connected,
so $M = M_{1} \oplus_{3} M_{2}$ for some matroids
$M_{1}$ and $M_{2}$.
Let $\{a,b,c\}$ be $E(M_{1})\cap E(M_{2})$.
Let $P$ be the generalised parallel connection of
$M_{1}$ and $M_{2}$, so that
$M=P\ba \{a,b,c\}$.
\Cref{prop10} implies that $M_{1}$ and
$M_{2}$ are isomorphic to proper minors of $M$.
Moreover $M_{1}$ and $M_{2}$ are loopless, and both
$\si(M_{1})$ and $\si(M_{2})$ are $3$\dash connected by
\cite[(4.3)]{Sey80}.
The following equality is from Walton and
Welsh~\cite{WW80}.
\begin{linenomath}
%\begin{multline}
\begin{equation}
\label{eqn8}
\chi(M; t) =
\frac{\chi(M_{1}; t)\chi(M_{2}; t)}{(t - 1)(t - 2)}
+ \chi(P \ba a \ba b / c; t)%\\
+ \chi(P \ba a / b; t) + \chi(P / a; t).
%\end{multline}
\end{equation}
\end{linenomath}
Since $P$ is binary, our earlier discussion means that
the evaluations
$\chi(P \ba a \ba b / c; 16)$,
$\chi(P \ba a / b; 16)$, and
$\chi(P / a; 16)$ are all non-negative.
On the other hand, by the minimality of $M$,
$\chi(M_{1};16)$ and $\chi(M_{2};16)$ are positive.
We deduce that $c(M;2)\leq 4$.
As $M$ is a counterexample, $c(M; 2) = 4$, so $\chi(M; 8) = 0$.
The terms in~\eqref{eqn8} must be zero when $t = 8$, so
we can assume by relabeling that $\chi(M_{1}; 8) = 0$.
Therefore $c(M_{1}; 2) = 4$.

The critical exponent of $\si(M_{1})$ is precisely the
critical exponent of $M_{1}$.
Since $\si(M_{1})$ is $3$\dash connected and obeys the
\namecref{thm19}, it follows that $\si(M_{1})\cong \pg{3,2}$.
Exactly as in the proof of \Cref{lem8}, we can show that $M$ has a
minor isomorphic to the matroid produced from \pg{3,2}
by performing a $\Delta\text{-}Y$ operation on $T$.
This matroid has an \mkt\dash minor
\cite[Appendix~C]{MRW10}, so we have a contradiction that
completes the proof.
\end{proof}

\section{Acknowledgements}

We thank the referees for their careful readings and very
constructive comments, one of which enabled a substantial
simplification.

%\bibliographystyle{/Users/Home/LaTeX/Bibliography/MRStyle}
%\bibliography{/Users/Home/LaTeX/Bibliography/References}

%\bibliographystyle{/Users/Dillon/Maths/LaTeX/Bibliography/MRStyle}
%\bibliography{/Users/Dillon/Maths/LaTeX/Bibliography/References}

\end{document}